\documentclass[12pt,a4paper]{article}
\usepackage{fullpage,amsfonts,amssymb,amsmath}

\author{Stephen J. Pride}
\date{}
\usepackage[all]{xy}

\def\Im{\mathop{\mathrm{Im}}}
\def\id{\mathop{\mathrm{id}}}
\begin{document}
\begin{center}
{\bf ON THE RESIDUAL FINITENESS AND OTHER PROPERTIES OF (RELATIVE)
ONE-RELATOR GROUPS}
\end{center}
%\maketitle
\begin{center}
{Stephen J. Pride}
\end{center}
{\textsc{Abstract.} A relative one-relator presentation has the
form {$\mathcal{P} = \langle \mathbf{x}, H; R \rangle$} where
$\mathbf{x}$ is a set, $H$ is a group, and $R$ is a word on
$\mathbf{x}^{\pm 1} \cup H$. We show that if the word on
$\mathbf{x}^{\pm 1}$ obtained from $R$ by deleting all the terms
from $H$ has what we call the {\it unique max-min property}, then
the group defined by $\mathcal{P}$ is residually finite if and
only if $H$ is residually finite (Theorem 1). We apply this to
obtain new results concerning the residual finiteness of
(ordinary) one-relator groups (Theorem 4). We also obtain results
concerning the conjugacy problem for one-relator groups (Theorem
5), and results concerning the relative asphericity of
presentations of the form $\mathcal{P}$ (Theorem 6).}\\

\noindent {\textit {2000 Mathematics Subject Classification.}
Primary 20E26, 20F05; Secondary 20F10, 57M07.\\

\noindent {\textit {Key words and phrases.} Residual finiteness,
one-relator group, relative presentation, (power) conjugacy
problem, asphericity, unique max-min property, 2-complex of
groups, covering complex.\\
\begin{center}
{\bf 1  Introduction}
\end{center}

The question of when one-relator groups are residually finite is
still open.

In the torsion-free case there are well-known examples of groups
which are not residually finite, namely the
Baumslag-Solitar/Meskin groups [4], [15]:

$$
G = \langle \mathbf{x}; U^{-1}V^{l}UV^{m} \rangle ,
$$
where $U$, $V$ do not generate a cyclic subgroup of the free group
on $\mathbf{x}$, and $|l| \neq |m|$, $|l|, |m| > 1$. On the other
hand, there are some examples which are known to be residually
finite. For instance, it was shown in [3] that if
$$
W = UV^{-1},\eqno{(1)}\
$$
where $U$, $V$ are positive words on an alphabet $\mathbf{x}$ and
the exponent sum of $x$ in $UV^{-1}$ is $0$ for each $x \in
\mathbf{x}$, or if

$$
W = [U, V],\eqno{(2)}\
$$
where $U$, $V$ are (not necessarily positive) words on
$\mathbf{x}$ such that no letter $x \in \mathbf{x}$ appears in
both $U$ and $V$, then $G = \langle \mathbf{x}; W \rangle$ is
residually finite.

In the torsion case there is the well-known open
question:\\

{\bf Question 1} [2], [5, Question OR1] Is every one-relator group
with torsion
residually finite?\\

Question $1$ is known to be true when $G = \langle \mathbf{x};
W^{n} \rangle$ where $W$ is a {\it positive} word and $n>1$ [9]
(see also [19]). In [20], Wise obtains further related results,
summed up by his ``Quasi-Theorem 1.3": {\it If $W$ is sufficiently
positive, and $W^n$ is sufficiently small cancellation, then $G$
is residually finite.}

A related open question is:\\

{\bf Question 2} [5,Question OR6], [11, Question 8.68] If a
torsion-free one-relator group $G_1 = \langle \mathbf{x}; W
\rangle$ is residually finite, then is $G_n = \langle \mathbf{x};
W^n \rangle$ also residually finite for
$n>1$?\\

(Of course, if Question 1 is true, then Question 2 is trivially
true.)

It was shown in [1] that Question 2 holds true when $W$ has the
form $(1)$ or $(2)$.

Here, amongst other things, we tackle Question 2 by considering
{\it relative} presentations.

A relative presentation has the form

$$
\mathcal{P} = \langle \mathbf{x}, H; \mathbf r \rangle
$$

\noindent where $H$ is a group and $\mathbf{r}$ is a set of
expressions of the form

$$
R = x_1^{\varepsilon_1} h_1 x_2^{\varepsilon_2} h_2 \ldots
x_r^{\varepsilon_r} h_r \  (r > 0, x_i \in
\mathbf{x},\varepsilon_i = \pm 1,\ h_i \in H , 1 \leq i \leq r).
\eqno{(3)}
$$

\noindent The word

$$
W = x_1^{\varepsilon_1} x_2^{\varepsilon_2} \ldots
x_r^{\varepsilon_r} \quad (r > 0, x_i \in \mathbf{x},
\varepsilon_i = \pm 1, 1 \leq i \leq r) \eqno{(4)}
$$

\noindent is called the $\mathbf{x}$-{\it{skeleton}} of $R$. We
{\it do not} require that the $\mathbf{x}$-skeleton is reduced or
cyclically reduced.
 The
group $G = G({\mathcal{P}})$ {\it defined by} $\mathcal{P}$ is the
quotient of $H*F$ (where $F$ is the free group on $\mathbf{x}$) by
the normal closure of the elements of $H*F$ represented by the
expressions $R \in\mathbf{r}$. The composition of the canonical
imbedding $H \rightarrow {H*F}$ with the quotient map $H*F
\rightarrow G$ is called the {\it natural homomorphism}, denoted
by $\nu : H \rightarrow G$ (or simply $H \rightarrow G$).

As is normal, we will often abuse notation and write $G = \langle
\mathbf{x}, H; \mathbf{r} \rangle$, or $G \cong \langle
\mathbf{x}, H; \mathbf{r} \rangle$.

When $\mathbf{r}$ consists of a single element $R$, then we have
the {\it one-relator relative presentation}
$$
\mathcal{P} = \langle \mathbf{x}, H; R \rangle. \eqno{(5)}
$$
\noindent Heuristically,  $G = G(\mathcal{P})$ should be governed
by the ``shape" of the $\mathbf{x}$-skeleton of $R$ and the
algebraic properties of $H$.

Here we introduce the {\it unique max-min property} for the
``shape" of $W$. (Words of the form (1) are a very special case.)
For a group $H$, denote by $\mathcal{M}_H$ the class of
one-relator relative presentations of the form (5), where
$W$ has the unique max-min property.\\

\noindent{\bf Theorem 1} {\it If $\mathcal{P}$ is in
$\mathcal{M}_H$ then}:

(i) {\it the natural homomorphism $H \rightarrow G(\mathcal{P})$
is injective};

(ii) $G(\mathcal{P})$ {\it is residually finite if and only if $H$ is residually finite}.\\

We can deduce from this\\

\noindent{\bf Theorem 2} (Substitution Theorem). {\it Let $K$ be a
one-relator group given by an ordinary presentation $\langle
\mathbf{y}, z; S(\mathbf{y},z)\rangle$, and let $\mathcal{P} =
\langle \mathbf{x}, H; R \rangle$ be an
$\mathcal{M}_H$-presentation. Then the group given by the relative
presentation $\langle \mathbf{x}, \mathbf{y}, H; S(\mathbf{y}, R)
\rangle$
is residually finite if and only if $H$ and $K$ are residually finite}.\\

We can give the proof of this straightaway.  Consider the
$\mathcal{M}_{H * K}$-presentation $\overline{\mathcal{P}}
=\langle \mathbf{x}, H *K; Rz^{-1} \rangle$. By Theorem 1, $L =
G(\overline{\mathcal{P}})$ is residually finite if and only if $H
* K$ is residually finite, which is equivalent to requiring that
both $H$ and $K$ are residually finite (using results discussed in
[12] p417). Now note that
$$
L \cong \langle \mathbf{x}, \mathbf{y}, z, H; S(\mathbf{y}, z),
Rz^{-1} \rangle \cong \langle \mathbf{x}, \mathbf{y}, H;
S(\mathbf{y}, R) \rangle.
$$

In particular, taking $K$ to be defined by  $\langle z; z^n \rangle \; (n > 1)$ we have:\\

\noindent{\bf Theorem 3}.  {\it If $G = \langle \mathbf{x}, H; R
\rangle$ is a residually finite $\mathcal{M}_H$-group, then the
group
$G_n = \langle \mathbf{x}, H; R^n \rangle \; (n > 1)$ is also residually finite}.\\

Now take $H$ to be a free group $\Phi$.  Then
$\mathcal{M}_\Phi$-groups are one-relator groups.  Since $\Phi$ is
residually finite ([12],p116 or p417), we obtain the following
theorem concerning residual finiteness of one-relator groups.\\

\noindent{\bf Theorem 4} {\it Every $\mathcal{M}_\Phi$-group $G =
\langle \mathbf{x}, \Phi; R \rangle$ is a residually finite
one-relator group.  Moreover, if $K = \langle \mathbf{y}, z;
S(\mathbf{y}, z)\rangle$ is a one-relator group, then the
one-relator group $\overline{K} = \langle \mathbf{x}, \mathbf{y},
\Phi; S(\mathbf{y}, R)\rangle$ is residually finite if and only if
$K$ is residually finite.  In particular,
$G_n = \langle \mathbf{x}, \Phi; R^n \rangle \; (n > 1)$ is residually finite}.\\

The solution of the conjugacy problem for one-relator groups with
{\it torsion} has been solved by B.B.Newman [16]. However, for the
{\it torsion-free} case the problem is still open [5, Question O5].\\

 \noindent{\bf Theorem 5} {\it Every
$\mathcal{M}_\Phi$-group ($\Phi$ a finitely generated free group)
has solvable conjugacy problem. Also, such groups have solvable {\it power} conjugacy problem.}\\

(Two elements $c, d$ of a group are said to be {\it power
conjugate} if some power of $c$ is conjugate to some power of
$d$.)
\\

Other aspects of relative presentations (and in particular,
one-relator relative presentations) have been studied intensively,
particularly {\it asphericity}. Recall [6] that a relative
presentation $\mathcal P$ is {\it aspherical} (more accurately,
{\it diagrammatically} aspherical) if every spherical picture over
$\mathcal P$ contains a dipole. Under a weaker condition on
``shape" (the {\it unique min property}, or equivalently the {\it
unique max property}) we can prove:\\

\noindent{\bf Theorem 6} {\it Let $\mathcal P$ be a relative
presentation as in (5), where $\it W$ has the unique min property.
Then $\mathcal P$ is aspherical.}\\

It then follows from [6] (see Corollary 1 of Theorem 1.1, Theorem
1.3, and Theorem 1.4) that for the group $G = G({\mathcal P})$ we
have:\\

(i) {\it the natural homomorphism $H \rightarrow G$ is
injective};\\

(ii) {\it every finite subgroup of $G$ is contained in a
conjugate of $H$};\\

(iii) {\it for any left $\mathbb{Z}G$-module $A$, and any right
$\mathbb{Z}G$-module $B$},
$$
H^n(G,A) \cong H^n(H,A),
$$
$$
H_n(G,B) \cong H_n(H,B)
$$

{\it for all $n \geq 3$}.\\

\begin{center}
{\bf 2  Max-min property}
\end{center}

Let $\mathbf{x}$ be an alphabet.  A {\it weight function} on
$\mathbf{x}$ is a function

$$
\theta: \mathbf{x} \longrightarrow \mathbb{Z}
$$

\noindent such that $\Im \theta$ generates the additive group
$\mathbb{Z}$ (that is, $gcd \{\theta(x): x \in \mathbf{x}\}$ is
1). A {\it strict} weight function is one for which $\theta(x)
\neq 0$ for all $x \in \mathbf{x}$.

Let $W$ be a  word on $\mathbf{x}$ as in $(4)$. Given a weight
function $\theta$, we then have the function

$$
\phi = \phi_W^{\theta}: \{0, 1, 2, \ldots, r\} \rightarrow
\mathbb{Z},
$$

$$
\phi(j) = \sum^j_{i=0} \varepsilon_i \theta(x_i)
$$

\noindent(where $\phi(0)=0$ since the empty sum is taken to be
$0$). We will say that the weight function is {\it admissible} for
$W$ if $\phi(r) = 0$.

For visual purposes, it is useful to extend $\phi$ to a piecewise
linear function $\phi: [0, r] \rightarrow \mathbb{R}$, so that the
graph of $\phi$ in the interval $[j-1, j]$ is the straight line
segment joining the points $(j-1, \phi(j-1)), \; (j, \phi(j))\; (0
< j \leq r)$. We will informally refer to this graph as ``the
graph of $W$" (with respect to $\theta$).

A word $W$ as in (4) will be said to have the {\it unique max-min
property} if for some admissible strict weight function $\theta$,
the graph of $W$ has a unique maximum and a unique minimum.  To be
precise, we require that, for some admissible strict weight
function, and some $k, l \in \{1, 2, \ldots, r\}$, we have
$\phi(j) < \phi(k)$ for all $j \in \{1,2, \ldots, r\}-k$ and
$\phi(j) > \phi(l)$ for all $j \in \{1, 2, \ldots, r\}-\{l\}$.  We
also require that $x_k \neq x_{k+1}$ and $x_l \neq x_{l+1}$
(subscripts modulo $r$). This amounts to requiring that $W$ is
``reduced at the unique maximum and minimum", that is,
$x_k^{\varepsilon_k} \neq x_{k+1}^{-\varepsilon_{k+1}},
x_l^{\varepsilon_l} \neq x_{l+1}^{-\varepsilon_{l+1}}$ (subscripts
modulo $r$).
%{\it It is important to note that this amounts to
%requiring that} $x_k \neq x_{k+1}$ and $x_l \neq x_{l+1}$.
For at the maximum and minimum we must have {\it either} $x_j \neq
x_{j+1}$, {\it or} $x_j = x_{j+1}$ and $\varepsilon_j = -
\varepsilon_{j+1}$ $(j = k,l)$.
%\;(j = k,l)$ then we would have to have
%$\varepsilon_j = \varepsilon_{j+1}$, which would be impossible at
%a maximum or minimum.
If the two letters occurring at the unique maximum are not
disjoint from the two letters occurring at the unique minimum
(i.e. $\{x_k, x_{k+1}\} \cap \{x_l, x_{l+1}\}$ is not empty), then
we will say that $W$ has the {\it strong} unique max-min property.

 A word {\it W} as in (4)
will be said to have the {\it unique min property} if for some
strict weight function $\theta$, the graph of {\it W} has a unique
minimum (but not necessarily a unique maximum). The {\it unique
max property} is defined similarly, but is not really of interest
because replacing $\theta$ by $-\theta$ will convert this property
to the unique min property.

We let $\mathcal{M}^{\mathbf{1}}_H$ (respectively
$\mathcal{S}^{\mathbf{1}}_H)$ denote the subclass of
$\mathcal{M}_H$ consisting of relative presentations of the form
(5) for which $W$ has the unique max-min property (respectively,
the strong unique max-min property) with respect to the weight
function

$$
\mathbf{1}: \mathbf{x} \longrightarrow \mathbb{Z} \quad x \mapsto
1\; (x \in \mathbf{x}).
$$
\\
\noindent{\bf Lemma 1} {\it Every $\mathcal{M}_H$-group can be
embedded into
an $\mathcal{M}^{\mathbf{1}}_H$-group}.\\

{\bf Proof.} Let $G$ = $\langle \mathbf{x},H; R \rangle$ with $R$
as in (3), and suppose $W = x_1^{\varepsilon_1}
x_2^{\varepsilon_1} \ldots x_2^{\varepsilon_r}$ has the unique
max-min property with respect to some strict weight function
$\theta: \mathbf{x} \rightarrow \mathbb{Z}$.  We can assume
$\theta(x) > 0$ for all $x$.  For if $\theta(x) < 0$ then we can
replace $x$ by $x^{-1}$.

Let
$$
\mathbf{y} = \{y: y \in \mathbf{x}, \theta(y) > 1\},
$$

\noindent and let

$$
\hat{\mathbf{x}} = (\mathbf{x} - \mathbf{y}) \cup \{y_1, y_2,
\ldots, y_{\theta(y)}: y \in \mathbf{y}\}.
$$

\noindent Let $\hat{G} = \langle \hat{\mathbf{x}}, H;
\hat{R}\rangle$, where $\hat{R}$ is obtained from $R$ by replacing
each occurrence of $y^{\pm 1}$ by $(y_1 y_2 \ldots
y_{\theta(y)})^{\pm 1} \; (y \in \mathbf{y})$.  It is easy to see
that the $\hat{\mathbf{x}}$-skeleton $\hat{W}$ of $\hat{R}$ has
the unique max-min property with respect to $\mathbf{1}:
\hat{\mathbf{x}} \rightarrow \mathbb{Z}$. (The graph of $\hat{W}$
is obtained from that of $W$ by ``stretching" along the horizontal
axis.) Moreover, $G$ is embedded into $\hat{G}$, for we have the
retraction $\rho$ with section $\mu$:

\[
\xymatrix{ \hat{G} \ar@<0.4ex>[r]^\rho & \ar@<0.4ex>[l]^\mu G }
\quad \rho \mu = \textstyle{\id_G}
\]

$\rho: x \mapsto x\;(x \in \mathbf{x} - \mathbf{y}), \; y_1
\mapsto y, y_i \mapsto 1 \; (y \in \mathbf{y}, 1 < i \leq
\theta(y)), \; h \mapsto h (h \in H)$,

$\mu: x \mapsto x \; (x \in \mathbf{x} - \mathbf{y}), y \mapsto
y_1 y_2 \ldots
y_{\theta(y)} \;(y \in \mathbf{y}), h \mapsto h \; (h \in H)$.\\

\noindent{\bf Lemma 2} {\it Every
$\mathcal{M}^{\mathbf{1}}_H$-group can be embedded
into an $\mathcal{S}^{\mathbf{1}}_H$-group.}\\

\noindent{\bf Proof.}  Let $G = \langle \mathbf{x}, H; R \rangle$,
where the $\mathbf{x}$-skeleton $W$ of $R$ has the unique max-min
property with respect to the constant function $\mathbf{1}:
\mathbf{x} \rightarrow \mathbb{Z}$. Suppose the letters occurring
at the unique maximum are $a, b$, and those occurring at the
unique minimum are $c,d$. We can assume that $\{a, b\} \cap
\{c,d\}$ is empty, otherwise there is nothing to prove.

Let $\mathbf{y} = \mathbf{x} - \{a, b, c, d\}$, and introduce a
new alphabet

$$
\hat{\mathbf{x}} = \{a, b, c, d, e\} \cup \{y_1, y_2: y \in
\mathbf{y}\}.
$$

\noindent Let $\hat{R}$ be obtained from $R$ as follows.  For each
$y \in \mathbf{y}$, replace all occurrences of $y^{\pm 1}$ by
$(y_1y_2)^{\pm 1}$, and replace all occurrences of $a^{\pm 1}$
(respectively, $b^{\pm 1}$, $c^{\pm 1}, d^{\pm 1}$) by $(ea)^{\pm
1}$ (respectively, $(be)^{\pm 1}, (ec)^{\pm 1}, (de)^{\pm 1}$).
Let $\hat{G} = \langle \hat{\mathbf{x}}, H; \hat{R} \rangle$, and
let $\hat{W}$ be the word obtained from $\hat{R}$ by deleting all
terms from $H$.  The graph of $\hat{W}$ under the weight function
$\mathbf{1}: \hat{\mathbf{x}} \rightarrow \mathbb{Z}$ is the graph
of $W$ magnified by a factor of 2, and $e$ occurs at the unique
maximum and the unique minimum.  Moreover, $G$ is embedded into
$\hat{G}$ for we have the retraction $\rho$ with section $\mu$:

\[
\xymatrix{ \hat{G} \ar@<0.4ex>[r]^\rho & \ar@<0.4ex>[l]^\mu G }
\quad \rho \mu = \textstyle{\id_G}
\]

$\rho: z \mapsto z\; (z \in \{a,b,c,d\}), e \mapsto 1, y_1 \mapsto
y, y_2 \mapsto 1\; (y \in \mathbf{y}), h \mapsto h\;(h \in H),$

$\mu: a \mapsto ea, b \mapsto be, c \mapsto ec, d \mapsto de, y
\mapsto y_1 y_2 \;(y \in \mathbf{y}), \;
h \mapsto h\;(h \in H)$.\\

\noindent{\bf Remark 1} Note that in both the above proofs we have
$\mu \nu = \hat{\nu}$, where $\nu :H \rightarrow G$, $\hat{\nu} :
H \rightarrow \hat{G}$ are the natural homomorphims. Thus if
$\hat{\nu}$ is injective then so
is $\nu$.\\

\noindent{\bf Remark 2} Note also from the proof of the above two
lemmas we get that every $\mathcal{M}_H$-group is a retract of an
$\mathcal{S}^{\mathbf{1}}_H$-group.\\

\noindent{\bf Remark 3} The referee has brought my attention to
the work of K.S.Brown [8], which is concerned with whether a
homomorphism $\chi$ from a one-relator group $B = \langle
\mathbf{x}; W\rangle$ ($|\mathbf{x}|\geq 2$, $W$ as in (4) and
cyclically reduced) onto $\mathbb{Z}$ has finitely generated
kernel. Such a homomorphism is induced by a weight function
$\theta$ which is admissible for $W$. However, since $\theta$ need
not be strict, it is necessary to interpret the max-min property
more widely. Thus the unique maximum could be a ``plateau'': ie,
for some $k \in \{1, 2, \ldots, r\}$ we could have
$\phi(k)=\phi(k+1)$ and $\phi(j) < \phi(k)$ for all $j \in \{1,2,
\ldots, r\}-\{k, k+1\}$ (subscripts modulo $r$). Similarly, the
unique minimum could be a ``reverse plateau". Then according to
Brown [8], as restated in Theorem 2.2 of [13], $\ker{\chi}$ is
finitely generated if and only if $|\mathbf{x}|=2$,
%say $\mathbf{x}= {\{x,y\}}$,
and $W$ has the unique max-min property in the above sense with
respect to the corresponding weight function. In our work we could
also allow non-strict weight functions. However, for the most part
this can be avoided. For example, if the unique maximum is a
plateau with $x_{k}\neq x_{k+2}$ then we could transform it to a
genuine maximum by deleting $x_{k+1}$ from $\textbf{x}$ and
replacing $H$ by $H \ast \langle x_{k+1} \rangle$. However, if the
unique maximum is a plateau with $x_{k}=x_{k+2}$ then some of our
arguments need to be modified, which we leave as an exercise for
the reader.\\

\begin{center}
{\bf 3 A construction}
\end{center}

By a {\it 2-complex of groups} we mean a connected graph of groups
(in the sense of Serre [18]) with trivial edge groups, together
with a set of closed paths, which we call {\it defining paths}.
(These are essentially the ``generalized complexes" defined in \S
1 of [10], where more details can be found. Note however, that in
[10] a ``2-cell" $c(\alpha)$ consists of {\it all} cyclic
permutations of $\alpha^{\pm 1}$ for each one of our defining
paths $\alpha$. We specifically {\it do not} add these extra
paths. This makes no significant difference.)

Let $\mathcal{P}$ be as in (5), and let $\theta$ be an admissible
weight function for $W$. There is then an induced epimorphism

$$
\psi: G \rightarrow \mathbb{Z} \quad x \mapsto \theta(x) \;(x \in
\mathbf{x}), h \mapsto 0 \; (h \in H).
$$

\noindent We can construct a 2-complex of groups

$$
\mathcal{\tilde{P}} = \langle\Gamma , H_n\;(n \in \mathbb{Z}); \
(n, R) \;(n \in \mathbb{Z})\rangle
$$
whose fundamental group is isomorphic to the kernel $K$ of $\psi$.
The underlying graph $\Gamma$ has vertex set $\mathbb{Z}$, edges
$(n, x^\varepsilon)\; (n \in \mathbb{Z}, x \in \mathbf{x},
\varepsilon = \pm 1)$, and initial, terminal and inversion
functions $\iota, \tau,^{-1}$ given by $\iota(n,x^\varepsilon)=n,
\tau(n,x^\varepsilon) = n+ \varepsilon\theta(x), (n,
x^\varepsilon)^{-1} = (n+ \varepsilon\theta(x), x^{-
\varepsilon})$. The vertex groups are copies $H_n = \{(n,h): h \in
H\}$ of $H$ (with the obvious multiplication $(n,h)(n,h') = (n,
hh'))$. We extend $\iota, \tau, ^{-1}$ to the elements of the
vertex groups by defining $\iota(n,h) = n = \tau(n,h), \;
(n,h)^{-1} = (n,h^{-1})$ (where $h^{-1}$ is the inverse of $h$ in
$H$).  We extend $\theta$ to $\mathbf{x}^{\pm 1} \cup H$ by
defining $\theta(x^{-1}) = - \theta(x)$  $(x \in \mathbf{x})$,
$\theta(h) = 0\; (h \in H)$. Then for any sequence $\alpha =
z_1z_2 \ldots z_q$ with $z_i \in \mathbf{x}^{\pm 1} \cup H$ and
any vertex $n \in \Gamma$, we have a path $(n, \alpha)$ in the
graph of groups starting at $n$, where
$$(n, \alpha) = (n, z_1)(n+\theta(z_1),z_2)(n+\theta(z_1)+\theta(z_2),z_3) \ldots
(n+\theta(z_1)+\theta(z_2) + \ldots + \theta(z_{q-1}), z_q).$$
\noindent In particular we have the (closed) paths $(n, R)$.

There is an obvious action of $\mathbb{Z}$ on the above graph of
groups, with $i \in \mathbb{Z}$ acting on vertices by $i \cdot n =
i+n\; (n \in \mathbb{Z})$, and on the edges and vertex groups by
$i.(n,z) =$\break $(i+n, z)\;(n \in \mathbb{Z}, z \in
\mathbf{x}^{\pm 1} \cup H)$.  This action of course extends to
paths.  Thus $(i, \alpha) = i.(0, \alpha)$. In particular, $(i,R)
= i.(0,R)$, so $\mathbb{Z}$ acts on $\mathcal{\tilde{P}}$.

If we regard $\mathcal{P}$ as a 2-complex of groups with a single
vertex $o$, edges $x^{\varepsilon}$ ($x \in \mathbf{x},
\varepsilon = \pm 1$), vertex group $H$, and defining path $R$,
then we have a mapping of 2-complexes of groups
$$
\rho: \mathcal{\tilde{P}} \longrightarrow \mathcal{P}
$$
$$
\ n \mapsto o,\ (n, x^\varepsilon)\; \mapsto x^\varepsilon,\ (n,
h)\; \mapsto h,\ (n, R)\mapsto R
$$
$(n \in \mathbb{Z}, x \in \mathbf{x}, \varepsilon = \pm 1, h \in
H)$. This induces a homomorphism

$$
\rho_\ast: \pi_1(\mathcal{\tilde{P}}, 0) \longrightarrow
\pi_1(\mathcal{P}, o) = G
$$
which is injective, and Im$\rho_\ast = K$. This can easily be
proved by adapting the standard arguments of covering space theory
for
ordinary 2-complexes (see for example [17] pp 157-159), to this relative situation.\\

\begin{center}
{\bf 4  Proof of Theorem 1}
\end{center}

Since residual finiteness is closed under taking subgroups, it
follows from Lemmas 1 and 2 and the Remark 1 at the end of \S 2
that it suffices to prove Theorem 1 for
$\mathcal{S}^{\mathbf{1}}_H$-groups.

We will make use of the following results: {\it (a) A free product
$F*B$, where $F$ is a free group, is residually finite if and only
if $B$ is residually finite;} {\it (b) An infinite cyclic
extension of a finitely generated group $L$ is residually finite
if and only if $L$ is residually finite.} (The first of these
follows from results on p417 of [12]; the second is a special case
of Theorem 7, p29 of [14].)

We can assume $\textbf{x}$ is finite. For if not let $\textbf{x}'$
be the set of letters occurring in $R$. Then $G$ is isomorphic to
$G'*\Psi$ where $G' \cong \langle \textbf{x}', H ; R\rangle$, and
$\Psi$ is the free group on $\textbf{x}-\textbf{x}'$. So by (a)
above, it is enough to work with $G'$.

Let $G$ be defined by an $\mathcal{S}^{\mathbf{1}}_H$
presentation as in (5), with $e \in \mathbf{x}$ occurring at both
the unique maximum and the unique minimum of the graph of $W$
under the weight function $\theta = \mathbf{1}$.  We denote the
maximum and minimum values of $\phi_W$ by $M$, $m$ respectively.
Note that $m \leq 0 \leq M$ and $m < M$.

We first deal with the trivial case when $M-m=1$. Then up to
cyclic permutation and inversion, $R=eha^{-1}h'$, where $a \in
\mathbf{x} - \{e\}$, $h,h' \in H$. Thus $G=\Phi*H$, where $\Phi$
is the free group on $\mathbf{x} - \{e\}$, so the theorem holds by
(a) above.

Now suppose $M-m>1$. Let $f \in \mathbf{x} - \{e\}$.

We have the epimorphism
$$
\psi: G \rightarrow \mathbb{Z} \quad x \mapsto 1 \;(x \in
\mathbf{x}), h \mapsto 0 \; (h \in H).
$$

\noindent Also, we have the homomorphism

$$
\eta: \mathbb{Z} \rightarrow G \quad 1 \mapsto f.
$$

\noindent Then $\psi \eta  = \id_{\mathbb{Z}}$, so $G$ is a
semidirect product $K \rtimes \mathbb{Z}$, where $K = \ker \psi$,
and with the action of $n \in \mathbb{Z}$ on $K$ being induced by
conjugation by $f^n$.

The fundamental group of $\mathcal{\tilde{P}}$ (at the vertex
$0$), as in \S 3, is isomorphic to $K$.

We will obtain a relative presentation for $K$ by collapsing a
maximal tree.

The edges $(n, f)^{\pm 1}$ form a maximal tree $T$ in $\Gamma$.
Let $R_n$ be the word on\break $\{(i, x): i \in \mathbb{Z}, x \in
\mathbf{x}, x \neq f\} \cup (\displaystyle{\bigcup_{i \in
\mathbb{Z}}}H_i)$ obtained from $(n, R)$ by deleting all edges
from $T$ which occur in $(n, R)$ and replacing all terms $(i,
x^{-1})$ by $(i-1, x)^{-1} \; (i \in \mathbb{Z}, x \in \mathbf{x},
x \neq f)$.  Then

$$
\mathcal{Q} = \langle (n, x)\; (n \in \mathbb{Z}, x \in
\mathbf{x}, x \neq f) \, , *_{n \in \mathbb{Z}}H_n \, ; R_n \;(n
\in \mathbb{Z})\rangle
$$

\noindent is a relative presentation for $K$.  Moreover, since the
edges in $T$ constitute an orbit under the action of $\mathbb{Z}$
on our graph of groups, the action of $\mathbb{Z}$ on $K$ is given
by the automorphism
$$
\mu: (n, x) \mapsto (n+1, x) \; (x \in \mathbf{x}, x \neq f), \;
(n, h) \mapsto (n+1, h)\; (h \in H)
$$
$(n \in \mathbb{Z})$.

Now consider the $HNN$-extension $\overline{K}$ of $K$ given by
the relative presentation

$$
\begin{array}{llll}
\overline{\mathcal{Q}} =&\langle(n, x)\; (n \in \mathbb{Z}, x \in
\mathbf{x},
x \neq f), *_{n \in \mathbb{Z}}H_n \, , s; R_n \; (n \in \mathbb{Z})\\[5pt]
& s(n,x)s^{-1} = (n+1, x) \; (n \in \mathbb{Z}, x \in \mathbf{x}, x \neq e, f),\\[5pt]
& s(n,h)s^{-1} = (n+1, h) \; (n \in \mathbb{Z}, h \in H) \rangle.
\end{array}
$$

\noindent The automorphism $\mu$ of $K$ can be extended to an
automorphism $\overline{\mu}$ of $\overline{K}$ by defining
$\overline{\mu}(s) = s$.  Then $G = K \rtimes_{\mu} \mathbb{Z}$
can be embedded into $\overline{G} = \overline{K}
\rtimes_{\overline{\mu}} \mathbb{Z}$.

By our assumption, up to cyclic permutation and inversion, $(0,
R)$ will have the form

$$
(M-1, e)(M,h)(M-1,a)^{-1}
\gamma_0((m,b)^{-1}(m,h')(m,e))^\varepsilon \delta_0,
$$

\noindent where $h, h' \in H, \varepsilon = \pm 1, a,b \in
\mathbf{x} - \{e\}$, and each term $(i,z)$ occurring in the paths
$\gamma_0, \delta_0$ is such that both its initial and terminal
vertices lie in the range $m+1, m+2, \ldots, M-1$.

Then
$$
R_0 = (M-1, e) \alpha_0 (m,e)^\varepsilon \beta_0
$$

\noindent where $\alpha_0, \beta_0$ do not contain any occurrence
of $(i,e)^{\pm 1}$ with $i \leq m$ or $i \geq M-1$. More
generally, for $n \in \mathbb{Z}$
$$
R_n = (n+M-1, e) \alpha_n (n+m,e)^\varepsilon \beta_n
$$
\noindent where $\alpha_n, \beta_n$ do not contain any occurrence
of $(i,e)^{\pm 1}$ with $i \leq n+m$ or $i \geq n+M-1$.

Let $F_0$ be the free group on
$$
(\mathbf{x} - \{e, f\}) \cup \{s, (m+1, e), (m+2, e) \ldots ,(M-1,
e)\}.
$$

\noindent Then there is a homomorphism
$$
\overline{K} \rightarrow H * F_0
$$
defined as follows:\\

%$s \mapsto s, (n,x) \mapsto s^n xs^{-n}\; (x \in \mathbf{x}, x
%\neq e, f, \; n \in \mathbb{Z})$,

%$(n,h) \mapsto s^nhs^{-n} \;(h \in H, n \in \mathbb{Z}), (i, e)
%$\mapsto (i,e) \quad
%(m+1 \leq i \leq M-1)$,\\
$$
\begin{array}{llll}
s \mapsto s, \\[5pt]
(n,x) \mapsto s^n xs^{-n}\; (x \in \mathbf{x}, x \neq e, f, \; n
\in \mathbb{Z}),\\[5pt]
(n,h) \mapsto s^nhs^{-n} \;(h \in H, n \in \mathbb{Z}),\\[5pt]
 (i, e)\mapsto (i,e) \quad
(m+1 \leq i \leq M-1),
\end{array}
$$
\noindent and (inductively), for $k=0,1,2, \ldots$

$$
(k+M,e) \mapsto \beta^{-1}_{k+1} (k+1+m, e)^{- \varepsilon}
\alpha_{k+1}^{-1},
$$

$$
(-k+m,e) \mapsto (\beta_{-k}(-k+M-1, e)\alpha_{-k})^{-
\varepsilon}.
$$

\noindent This homomorphism is actually an isomorphism.  The
inverse is defined by

$$
\begin{array}{llll}
x &\mapsto& (0, x) \quad (x \in \mathbf{x}, x \neq e, f),\\[5pt]
h &\mapsto& (0, h) \quad (h \in H),\\[5pt]
(i,e) &\mapsto& (i, e) \quad m+1 \leq i \leq M-1,\\[5pt]
s &\mapsto& s.
\end{array}
$$

Thus $\overline{G}$ is an infinite cyclic extension of the group
$F_{0} * H$.\\

\noindent{\bf Remark 4} Note that by sending $s$ to the generator
$1 \in \mathbb{Z} \subset {G = K \rtimes_{\mu} \mathbb{Z}}$, we
obtain a retraction of $\overline {G}$ onto $G$ (with section
induced by the inclusion of $K$ into $\overline {K}$).\\

We can now complete the proof.

Clearly the natural homomorphism from $H$ into $\overline{G}$ is
injective (and is thus injective into $G$). Hence if $H$ is not
residually finite then neither is $G$. It remains to show that if
$H$ is residually finite then so is $\overline{G}$ (and thus $G$).

\textit{Case 1.} If $H$ is finitely generated then the result
holds straight away by (a) and (b) above.

\textit{Case 2.} Suppose that $H$ is not finitely generated. For
any homomorphism $\theta$ from $H$ to a group $H_{\theta}$ we
obtain an induced homomorphism from $\overline{G}= (F_{0} *
H)\rtimes_{\overline{\mu}} \mathbb{Z}$ to $\overline{G}_{\theta}=
(F_{0} * H_{\theta})\rtimes_{\overline{\mu}} \mathbb{Z}$ which
acts as $\theta$ on $H$ and acts as the identity on $F_0$ and
$\mathbb{Z}$. Let $g=(w_{0}h_{1}w_{1} \ldots h_{q}w_{q}).n$  be a
non-trivial element of $\overline{G}$ (where $q \geq 0, h_{1}
\ldots h_{q} \in H - {\{1}\}, w_{1}, \ldots ,w_{q-1} \in F_{0} -
{\{1}\}, w_{0}, w_{q} \in F_{0}$, $n \in \mathbb{Z}$, and if $q$
is $0$ then either $n \neq 0$ or $w_{0}$ is non-trivial). Since
residually finite groups are fully residually finite, there is a
homomorphism $\tau$ from $H$ onto a finite group $H_{\tau}$ such
that $\tau (h_i)\neq 1$ $(i=1,\ldots ,q)$. So the image of $g$ in
$\overline{G}_{\tau}= (F_{0} * H_{\tau})\rtimes_{\overline{\mu}}
\mathbb{Z}$ is non-trivial, and then Case 1 applies.
%$\overline{G}$ is
%residually finite if
%and only if $H$ is residually finite.

\begin{center}
{\bf 5 Proof of Theorem 5}
\end{center}

\noindent{\bf Lemma 3} {\it Let $C$ be a group which is a retract
of a group $B$. If $B$ has solvable conjugacy (or power conjugacy)
problem, then so does $C$.}

{\bf Proof.} By assumption we have maps $\xymatrix{ {B}
\ar@<0.4ex>[r]^\rho & \ar@<0.4ex>[l]^\mu C }, \rho \mu =
\textstyle{\id_C}$. Clearly if $c, d \in {C}$ are conjugate
(respectively, power conjugate) in $C$ then $\mu{(c)}, \mu{(d)}$
are conjugate (respectively, power conjugate) in $B$. Conversely
if there exists $b \in B$ such that $b \mu {(c)} {b^{-1}}=\mu
{(d)}$ (respectively, $b {\mu {(c)}^{i}} {b^{-1}}={\mu
{(d)}^{j}}$), then ${\rho{(b)}} {c} {\rho{(b)}^{-1}} = {d}$
(respectively, ${\rho{(b)}} {c^i} {\rho{(b)}^{-1}} = {d^j}$). Thus
the result
follows.\\

Now it is shown in [7] that infinite cyclic extensions of finitely
generated free groups have solvable conjugacy, and power
conjugacy, problem. By Remarks  2, 4, every
$\mathcal{M}_\Phi$-group is a retract of such a group.\\

\begin{center}
{\bf 6 Proof of Theorem 6}
\end{center}

We will assume familiarity with the terminology in \S\S 1.2, 1.4
of [6].

As in Lemma 1, we can assume that $\theta(x)> 0$ for all $x$. We
can extend $\theta$ to any word $U = y_1^{\varepsilon_1}
y_2^{\varepsilon_2} \ldots y_s^{\varepsilon_s}$, $(s>0, y_i \in
{\mathbf{x}}, \varepsilon_{i} =\pm 1, 1 \leq i \leq s)$ by $\theta
(U) = \sum^s_{i=0} \varepsilon_i \theta(y_i)$.

Let $\mathbb{P}$ be a based connected spherical picture (with at
least one disc) over $\mathcal{P}$, with global basepoint $O$, and
basepoint $O_{\Delta}$ for each disc $\Delta$. (Note that since
$R$ is not periodic, there will be just one basepoint for each
disc.) We will also choose, for each region $\textsf{R}$, a point
$O_{\textsf{R}}$ in the interior of $\textsf{R}$.

We can relabel $\mathbb P$ to obtain a picture
$\mathbb{\tilde{P}}$ over $\mathcal{\tilde{P}}$ as follows:

(a) For each region $\textsf{R}$, choose a tranverse path
$\gamma_\textsf{R}$ from $O$ to $O_ \textsf{R}$, and let
$U_\textsf{R}$ (a word on $\mathbf{x}$) be the label on the path
$\gamma_\textsf{R}$. Then the {\it potential} $q(\textsf{R})$ of
$\textsf{R}$ is $\theta({U_\textsf{R}})$. (This is independent of
the choice of path $\gamma_\textsf{R}$, since $\theta(W) = 0$.)

(b) For an arc tranversely labelled $x \in \mathbf{x}$ say,
relabel it by $(q(\textsf{R}), x)$ where $\textsf{R}$ is the
region where the tranverse arrow on the arc begins.

(c) For a corner of a disc, with label $h \in H$ say, relabel the
corner by $(q, h)$, where $q$ is the potential of the region in
which the corner occurs.

For a disc $\Delta$, let $q_{\Delta}$ be the potential of the
region containing $O_{\Delta}$. Then in the relabelled picture,
$\Delta$ will be labelled by the path $(q_{\Delta}, R)$.

Let $\Theta$ be a {\it minimal} disc, that is, a disc such that
$q_{\Theta}\leq q_{\Delta}$ for all discs $\Delta$. Let $m$ be the
minimum value of $\phi_W^{\theta}$, and let $e$ be one of the two
distinct letters occurring at the unique minimum. Then in the path
$(0, R)$ there is a unique edge labelled $(m, e)$. Now $\Theta$ is
labelled by $(q_{\Theta}, R)$ in $\tilde{\mathbb{P}}$, and thus
there is a unique edge labelled $(m+q_{\Theta}, e)$ incident with
$\Theta$. This arc must intersect another disc $\Theta^{'}$, which
must also be labelled by $(q_{\Theta}, R)$, but with the opposite
orientation. Thus we obtain a dipole in $\mathbb{\tilde{P}}$ where
$\Theta$, $\Theta^{'}$ are the discs of the dipole. Reverting to
$\mathbb{P}$, this dipole in $\mathbb{\tilde{P}}$ gives rise to a
dipole in $\mathbb{P}$.\\

\noindent \textit{Acknowledgement}. I thank the referee for
his/her helpful comments.

\vspace{0.50cm}
\begin{center}
{\bf References }
\end{center}
\noindent [1] R.B.J.T.Allenby and C.Y.Tang, Residual finiteness of
certain 1-relator groups: extensions of results of Gilbert
Baumslag, {\it Math. Proc. Camb. Phil. Soc.} $\mathbf{97}$ (1985),
225-230.

\noindent [2] G.Baumslag, Residually finite one-relator groups,
{\it Bull. Amer. Math. Soc.} $\mathbf{73}$ (1967), 618-620.

\noindent [3] G.Baumslag, Free subgroups of certain one-relator
groups defined by positive words, {\it Math. Proc. Camb. Phil.
Soc.} $\mathbf{93}$ (1985), 247-251.

\noindent [4] G.Baumslag and D.Solitar, Some two-generator
one-relator non-Hopfian groups, {\it Bull. Amer. Math. Soc.}
$\mathbf{68}$ (1962), 199-201.

\noindent [5] G.Baumslag, A.Miasnikov and V.Shpilrain, Open
problems in combinatorial and geometric group theory,
http://zebra.sci.ccny.cuny.edu/web/nygtc/problems/

\noindent [6] W.A.Bogley and S.J.Pride, Aspherical relative
presentations, {\it Proc. Edin. Math. Soc.} $\mathbf{35}$ (1992),
1-39.

\noindent [7] O.Bogopolski, A.Martino, O.Maslakova and E.Ventura,
The conjugacy problem is solvable for free-by-cyclic groups, {\it
Bull. London Math. Soc.} $\mathbf{38}$ (2006), 787-794.

\noindent [8]  K.S.Brown, Trees, valuations, and the
Bieri-Neumann-Strebel invariant, {\it Invent. Math.}
$\textbf{90}$, (1987), 479-504.

\noindent [9] V.Egorov, The residual finiteness of certain
one-relator groups, {\it Algebraic Systems, Ivanov. Gos.Univ.,
Ivanovo} (1981), 100-121.

\noindent [10] J.Howie and S.J.Pride, A spelling theorem for
staggered generalized 2-complexes, with applications, {\it Invent.
Math.} {\bf 76} (1984), 55-74.

\noindent [11] Kourovka Notebook $\mathbf{15}$ (2002).

\noindent [12] W.Magnus, A.Karrass and D.Solitar, {\it
Combinatorial Group Theory} (Second Edition), Dover, New York,
1976.

\noindent [13] J.Meier, Geometric invariants for Artin groups,
{\it Proc. London Math. Soc.(3)} {\bf 74} (1997), 151-173.

\noindent [14] C.F.Miller III, {\it On group-theoretic decision
problems and their classification}, Annals Of Mathematics Studies
{\bf 68}, Princeton University Press, 1971.

\noindent [15] S.Meskin, Nonresidually finite one-relator groups,
{\it Trans. Amer. Math. Soc.} $\mathbf{164}$ (1972), 105-114.

\noindent [16] B.B. Newman, Some results on one-relator groups,
{\it Bull. Amer. Math. Soc.} $\mathbf{74}$ (1968), 568-571.

\noindent [17] S.J.Pride, Star-complexes, and the dependence
problems for hyperbolic complexes, {Glasgow Math. J.} {\bf 30}
(1988), 155-170.

\noindent [18] J.-P.Serre, {\it Trees}, Springer-Verlag, Berlin
Heidelberg New York, 1980.

\noindent [19] D.Wise, The residual finiteness of positive
one-relator groups, {\it Comment. Math. Helv.} $\mathbf{76}$
(2001), 314-338.

\noindent [20] D.Wise, Residual finiteness of quasi-positive
one-relator groups, {\it J. London Math. Soc.}(2) $\mathbf{66}$
(2002), 334-350.\\

\vspace{10mm}

\noindent {\it Address}:

\noindent Department of Mathematics, University of Glasgow,
University Gardens, Glasgow G12 8QW, Scotland UK\\

\noindent sjp@maths.gla.ac.uk

\noindent
\end{document}